\documentclass[reqno,11pt]{amsproc}

\usepackage[dvips]{graphicx}
\usepackage{amsmath,amssymb}
\begin{document}

\title[On the construction of de Branges spaces]{On the construction of de Branges spaces for
dynamical systems associated with finite Jacobi matrices.}

\author[A.\,S.~Mikhaylov, V.\,S.~Mikhaylov]
{$^1$, $^2$A.\,S.~Mikhaylov, $^1$,V.\,S.~Mikhaylov}

\address{
$^1$ St. Petersburg Department of V.A. Steklov Institute of
Mathematics of the Russian Academy of Sciences, 7, Fontanka,
191023 St. Petersburg, Russia. $^2$ Saint Petersburg State
University, St.Petersburg State University, 7/9 Universitetskaya
nab., St. Petersburg, 199034 Russia.}

\email{mikhaylov@pdmi.ras.ru, vsmikhaylov@pdmi.ras.ru}

\begin{abstract}
We consider dynamical systems with boundary control associated
with finite Jacobi matrices. Using the method previously developed
by the authors, we associate with these systems special Hilbert
spaces of analytic functions (de Branges spaces).
\end{abstract}

\keywords{Boundary control method, Krein equations, Jacobi
matrices, De Branges spaces}

\maketitle

\newtheorem{corollary}{Corollary}
\newtheorem{definition}{Definition}
\newtheorem{lemma}{Lemma}
\newtheorem{proposition}{Proposition}
\newtheorem{remark}{Remark}
\newtheorem{theorem}{Theorem}
%\numberwithin{equation}{section}
%\newcommand{\thmref}[1]{Theorem~\ref{#1}}
%\newcommand{\secref}[1]{\S\ref{#1}}
%\newcommand{\lemref}[1]{Lemma~\ref{#1}}

\section{Introduction}

For a given sequence of positive numbers $\{a_0,$
$a_1,\ldots,a_{N-1}\}$ (in what follows we assume $a_0=1$) and
real numbers $\{b_1, b_2,\ldots,b_N \}$, we denote by $A$ the
finite Jacobi matrix given by
\begin{equation}
\label{Jac_matr}
A=\begin{pmatrix} b_1 & a_1 & 0 & 0 & 0 &\ldots \\
a_1 & b_2 & a_2 & 0 & 0 &\ldots \\
\ldots &\ldots  &\ldots &\ldots & \ldots &\ldots \\
\ldots &\ldots  &\ldots &a_{N-2} & b_{N-1} &a_{N-1}\\
\ldots &\ldots &\ldots &\ldots & a_{N-1} &b_N
\end{pmatrix}.
\end{equation}
Let $u=(u_1,\ldots,u_N)\in \mathbb{R}^N$ and $T>0$ be fixed. With
the matrix $A$ we associate the dynamical system:
\begin{equation}
\label{DynSyst}
\begin{cases}
u_{tt}(t)-Au(t)=F(t),\quad t>0,\\
u(0)=u_t(0)=0,
\end{cases}
\end{equation}
where the vector function $F(t)=(f(t),0,\ldots,0)$, $f\in
L_2(0,T)$ is interpreted as a \emph{boundary control}. The
solution of (\ref{DynSyst}) is denoted by $u^f$. With the system
(\ref{DynSyst}) we associate the \emph{response operator} acting
by the rule
\begin{equation}
\label{RespJM} \left(R^Tf\right)(t)=u^f_1(t),\quad 0<t<T.
\end{equation}
The forward and inverse problems for the system (\ref{DynSyst})
and for the special case of this system, the finite
Krein-Stieltjes string, were the subjects of \cite{MM3,MM4}, where
as a main tool we used the Boundary control method \cite{B07,B08}.
In this paper we would like to demonstrate one more application of
the Boundary control method, namely the construction of the de
Banges space associated with (\ref{DynSyst}).

De Branges spaces play an important role in the inverse spectral
theory of first order canonical systems, see for example
\cite{DBr,DMcK,RR}. In \cite{MM6,MM7} the authors shows how to use
the Boundary control method to associate de Branges spaces with
different dynamical systems. Note that our approach differs from
the classical one and potentially admits the generalization to
multidimensional systems. The algorithm proposed in \cite{MM7,MM6}
is as follows: fixing some finite time $t=T$ one denote by
$\mathcal{F}^T$ the set of controls acting on the time interval
$(0,T)$ and introduces the \emph{reachable set} of the dynamical
system at this time:
\begin{equation*}
U^T:=\{u^f(T)\,|\,f\in \mathcal{F}^T\}.
\end{equation*}
Then one applies the Fourier transform $\mathbb{F}$ associated
with the operator $A$ to elements from $U^T$ and get a linear
manifold $\mathbb{F}U^T$. Then this linear manifold is equipped
with the norm defined by the \emph{connecting operator}, which
resulted in the de Branges space associated with the initial
dynamical system. In the considered in \cite{MM6,MM7} models the
dynamical systems have different properties in respect to the
\emph{boundary controllability}: the system associated with the
Schr\"odinger operator is exactly controllable from the boundary,
the system associated with the one-dimensional Dirac operator is
not controllable, but the controllability restores after some
trick associated with doubling the state of the system, the
discrete system associated with a finite Jacobi matrix is boundary
controllable, but the time in the model considered was discrete,
see also \cite{MM11} for the case of semi-infinite matrix. The
peculiarity of the system (\ref{DynSyst}) is the lack of the
boundary controllability, and in opposite to all systems
considered in \cite{MM7,MM11}, the speed of the wave propagation
in (\ref{DynSyst}) is infinite. Nevertheless we will show that the
method from \cite{MM6,MM6} works, and using it one can construct
de Branges space associated with $A$.

In the second section we provide the necessary information on the
solution of the forward and inverse problems for (\ref{DynSyst})
from \cite{MM4}. In the third section we remind the reader some
useful definitions and construct the de Branges space associated
with (\ref{DynSyst}).

\section{Dynamical system, forward problem, Krein equations.}

The following Cauchy problem for the difference equation
\begin{equation}
\label{Cauchy}
\begin{cases}
a_1\phi_2+b_1\phi_1=\lambda\phi_1,\\
a_n\phi_{n+1}+a_{n-1}\phi_{n-1}+b_n\phi_n=\lambda \phi_n,\quad
n=2,\ldots,N,\\
\phi_1=1,
\end{cases}
\end{equation}
determines the set of polynomials
$\{1,\phi_2(\lambda),\ldots,\phi_N(\lambda),\phi_{N+1}(\lambda)\}$.
Let $\lambda_1,\ldots,\lambda_{N}$ be the roots of the equation
$\phi_{N+1}(\lambda)=0$, it is known \cite{AH} that they are real
and distinct. We denote by $(\cdot,\cdot)$ the scalar product in
$\mathbb{R}^N$ and introduce the vectors and the coefficients by
the rules:
\begin{equation*}
\varphi(\lambda)=\begin{pmatrix} \phi_1(\lambda)\\
\cdot\\ \cdot\\ \phi_{N}(\lambda)\end{pmatrix},\quad\varphi_k=\begin{pmatrix} \phi_1(\lambda_k)\\
\cdot\\ \cdot\\ \phi_{N}(\lambda_k)\end{pmatrix},\quad
\rho_k=\left(\varphi_k,\varphi_k\right),\quad k=1,\ldots,N.
\end{equation*}
Thus $\varphi_k$ are non-normalized eigenvectors of $A$,
corresponding to eigenvalues $\lambda_k$:
\begin{equation*}
A\varphi_k=\lambda_k\varphi_k,\quad \quad k=1,\ldots,N.
\end{equation*}
We call by \emph{spectral data} and \emph{the spectral function}
$\rho$ the following objects:
\begin{equation*}
\label{measure_JM} \left\{\lambda_i,\rho_i\right\}_{i=1}^{N},\quad
\rho(\lambda)=\sum_{\{k:\lambda_k<\lambda\}}\frac{1}{\rho_k}.
\end{equation*}
%Note the property:
%\begin{equation}
%\label{rho_pr} \sum_{k=1}^N\frac{1}{\rho_k}=1.
%\end{equation}

The standard application of the Fourier method yields
\begin{lemma}
The solution to (\ref{DynSyst}) admits the spectral representation
\begin{equation}
\label{Sol_spec_repr_JM}
u^f(t)=\sum_{k=1}^{N}h_k(t)\varphi_k,\quad
u^f(t)=\int_{-\infty}^\infty\int_0^t
S(t-\tau,\lambda)f(\tau)\,d\tau\varphi(\lambda)\,d\rho(\lambda),
\end{equation}
where
\begin{eqnarray*}
h_k(t)=\frac{1}{\rho_k}\int_0^t
f(\tau)S_k(t-\tau)\,d\tau,\\
S(t,\lambda)=\begin{cases}
\frac{\sin{\sqrt{\lambda}t}}{\sqrt{\lambda}},\quad
\lambda>0,\\
\frac{\operatorname{sh}{\sqrt{|\lambda|}t}}{\sqrt{|\lambda|}},\quad
\lambda<0,\\
t,\quad \lambda=0,
\end{cases},\quad S_k(t)=S(t,\lambda_k).
\end{eqnarray*}
\end{lemma}

We introduce the \emph{outer space} of the system (\ref{DynSyst}),
the space of controls: $\mathcal{F}^T:=L_2(0,T;\mathbb{C})$ with
the scalar product $f,g\in \mathcal{F}^T,$
$\left(f,g\right)_{\mathcal{F}^T}=\int_0^Tf(t)\overline{g(t)}\,dt$.
The \emph{response operator} $R^T: \mathcal{F}^T\mapsto
\mathcal{F}^T$ is introduced by the formula (\ref{RespJM}). Making
use of (\ref{Sol_spec_repr_JM}) implies the representation formula
for $R^T$:
\begin{equation*}
%\label{Resp_def}
\left(R^Tf\right)(t)=u_1^f(t)=\sum_{k=1}^{N}h_k(t)=\int_0^t
r(t-s)f(s)\,ds,
\end{equation*}
where
\begin{equation*}
%\label{resp_func_JM}
r(t)=\sum_{k=1}^{N}\frac{1}{\rho_k}S_k(t),%\quad \sum_{k=1}^N\frac{1}{\rho_k}=1.
\end{equation*}
is called a \emph{response function}. Note that the operator $R^T$
is a natural analog of a dynamic Dirichlet-to-Neumann operator
\cite{B07} in continuous, and \cite{MM7,MM8,MM10} in discrete
cases.

The \emph{inner space} of (\ref{DynSyst}), i.e. the space of
states is $\mathcal{H}^N:=\mathbb{C}^{N}$, indeed for any $T>0$
and $f\in \mathcal{F}^T$, we have that $u^f(T)\in \mathcal{H}^N$.
The scalar product in $\mathcal{H}^N$ is given by
\begin{equation*}
(a,b)_{\mathcal{H}^N}=\sum_{k=1}^N a_k\overline{b_k}.
\end{equation*}
The \emph{control operator} $W^T: \mathcal{F}^T\mapsto
\mathcal{H}^N$ is introduced by the rule:
\begin{equation*}
W^Tf=u^f(T).
\end{equation*}
Due to (\ref{Sol_spec_repr_JM}) we have that
$W^Tf=\sum_{k=1}^{N}h_k(T)\varphi_k$. In \cite{MM3,MM4} the
authors used real inner and outer spaces, but in the complex case
all the results are valid as well.

%The boundary controllability properties of a dynamical system
%plays a crucial role in a procedure of solving the dynamic inverse
%problems, see \cite{B17,MM2,BM_1,AM,B08,MM1}. We set up the
%following \emph{control problem}: for a fixed state $a\in
%\mathcal{H}^N$, $a=\sum_{k=1}^{N}a_k\varphi_k$ we look for a
%control $f\in \mathcal{F}^T$ driving the system (\ref{DnSst}) to a
%prescribed state:
%\begin{equation}
%\label{control_JM} W^Tf=a.
%\end{equation}
%Using (\ref{Sol_spec_repr_JM}) we see that the equality
%(\ref{control_JM}) is equavalent to the following moment problem:
%to find $f\in \mathcal{F}^T$ such that
%\begin{equation*}
%a_k=h_k(T)=\frac{1}{\rho_k}\int_0^T f(\tau)S_k(T-\tau)d\tau,\quad
%k=1,\ldots,N.
%\end{equation*}
%Clearly, a solution to this moment problem exists, but it is not
%unique.
We introduce the subspace
\begin{equation*}
\mathcal{F}^T_1={\operatorname{Lin}\left\{S_k(T-t)\right\}_{k=1}^N},
\end{equation*}
where we assume complex coefficients in the span. The following
lemma establishes the boundary controllability of (\ref{DynSyst}):
\begin{lemma}
\label{LemmaCont_JM} The operator ${W}^T$ maps $\mathcal{F}^T_1$
onto $\mathcal{H}^N$ isomorphically.
\end{lemma}
%\begin{proof}
%Let $\left\{\widetilde S_k(T-t)\right\}_{k=1}^N,$ where
%$\widetilde S_k(T-t)\in \mathcal{F}^T_1$ be a bi-orthogonal basis
%in $\mathcal{F}^T_1$, i.e. $\int_0^T S_k(T-t)\widetilde
%S_l(T-t)\,dt=\delta_{kl}.$ Then we immediately have that
%\begin{equation*}
%f(t)=\sum_{k=1}^N a_k\rho_k\widetilde S_k(T-t).
%\end{equation*}
%\end{proof}

The \emph{connecting operator} $C^T:\mathcal{F}^T\mapsto
\mathcal{F}^T$ is defined by the rule
$C^T:=\left(W^T\right)^*W^T$, so by the definition for $f,g\in
\mathcal{F}^T$ one has
\begin{equation}
\label{C_T_JM}
\left(C^Tf,g\right)_{\mathcal{F}^T}=\left(u^f(T),u^g(T)\right)_{\mathcal{H}^N}=\left(W^Tf,W^Tg\right)_{\mathcal{H}^N}.
\end{equation}
It is crucial in the Boundary control method that $C^T$ can be
expressed in terms of inverse data:
\begin{theorem}
The connecting operator admits the representation in terms of
dynamic inverse data:
\begin{equation*}
%\label{CT_dyn_repr_JM}
\left(C^Tf\right)(t)=\frac{1}{2}\int_0^T\int_{|t-s|}^{2T-s-t}r(\tau)\,d\tau
f(s)\,ds,
\end{equation*}
and in terms of spectral inverse data:
\begin{equation}
\label{CT_sp_repr_JM} \left(C^Tf\right)(t)=\int_0^T
\sum_{k=1}^{N}\frac{1}{\rho_k}S_k(T-t)S_k(T-s)f(s)\,ds.
%C(t,s)=\frac{1}{l_1}\sum_{k=1}^{N-1}\frac{\sin{\sqrt{|\lambda_k|}(T-t)}\sin{\sqrt{|\lambda_k|}(T-s)}}{{|\lambda_k|}\rho_k}.\notag
%=\frac{1}{l_1}\int_{-\infty}^\infty
%\frac{\sin{\sqrt{|\lambda|}(T-t)}\sin{\sqrt{|\lambda|}(T-s)}}{{|\lambda|}}\,d\rho(\lambda)\notag
\end{equation}
\end{theorem}
%\begin{proof} Taking arbitrary $f,g\in \mathcal{F}^T$ and
%introducing the function
%$\psi(t,s):=\left(u^f(t),u^g(s)\right)_{\mathcal{H}^N}$, one can
%show that $\psi$ satisfies the equation
%\begin{equation*}
%\psi_{tt}-\psi_{ss}=\frac{1}{l_1}\left(f(t)(Rg)(t)-g(s)(Rf)(t)\right).
%\end{equation*}
%Solving it by the d'Alembert method and noting that
%$\psi(T,T):=\left(C^Tf,g\right)_{\mathcal{F}^T}$, yields
%(\ref{CT_dyn_repr_JM}). The formula (\ref{CT_sp_repr_JM}) follows
%from (\ref{C_T_JM}) and (\ref{Sol_spec_repr_JM}).
%\end{proof}

\begin{remark}
The formula (\ref{CT_sp_repr_JM}) implies that
$\mathcal{F}^T_1=C^T\mathcal{F}^T$, so $\mathcal{F}^T_1$ is
completely determined by inverse data. %The important properties of
%functions from $\mathcal{F}^T_1$ is that
%\begin{equation}
%\label{FT_proper_JM} f(T)=0,\quad \text{for all}\,\, f\in
%\mathcal{F}^T_1.
%\end{equation}
\end{remark}

\subsection{Krein equations.}
\label{Kr_eq}

By $f_k^T\in \mathcal{F}^T_1$ we denote the controls, driving the
system (\ref{DynSyst}) to prescribed \emph{special states}
\begin{equation*}
d_k\in \mathcal{H}^N, \,d_k=\left(0,\ldots,1,\ldots,0\right),\quad
k=1,\ldots,N.
\end{equation*}
%In other words, $f^T_k$ are the solutions to the equations
%$W^Tf_k^T=d_k$, $k=1\,\ldots,N$. By Lemma \ref{LemmaCont_JM} we
%know that such controls exist and are unique.
It is important that such a controls can be found as the solutions
to the Krein equations:
\begin{theorem}
\label{Krein_Th} The control $f_1^T$ can be found as the solution
to the following equation:
\begin{equation}
\label{Cont_f1_JM} \left(C^Tf_1^T\right)(t)=r(T-t),\quad 0<t<T.
\end{equation}
The controls $f^T_k$,  $k=2,\ldots,N$ satisfy the system:
\begin{equation}
\label{Cont_system_JM}
\begin{cases}
-\left(C^Tf_1^T\right)''=b_1C^Tf^T_1+a_1C^Tf^T_{2},\\
-\left(C^Tf_k^T\right)''=a_{k-1}C^Tf^T_{k-1}+b_kC^Tf^T_k+a_kC^Tf^T_{k+1},\quad
k=2,\ldots,N-1,\\
-\left(C^Tf_N^T\right)''=a_{N-1}C^Tf^T_{N-1}+b_NC^Tf^T_N.
\end{cases}
\end{equation}
%where we set $a_0=a_{N+1}=0.$
\end{theorem}

\section{De Branges space for $A$.}

Here we provide the information on de Branges spaces in accordance
with \cite{DBr,RR}. The entire function $E:\mathbb{C}\mapsto
\mathbb{C}$ is called a \emph{Hermite-Biehler function} if
$|E(z)|>|E(\overline z)|$ for $z\in \mathbb{C}_+$. We use the
notation $F^\#(z)=\overline{F(\overline{z})}$. The \emph{Hardy
space} $H_2$ is defined by: $f\in H_2$ if $f$ is holomorphic in
$\mathbb{C}^+$ and
$\sup_{y>0}\int_{-\infty}^\infty|f(x+iy)|^2\,dx<\infty$. Then the
\emph{de Branges space} $B(E)$ consists of entire functions such
that:
\begin{equation*}
B(E):=\left\{F:\mathbb{C}\mapsto \mathbb{C},\,F \text{ entire},
%\int_{\mathbb{R}}\left|\frac{F(\lambda)}{E(\lambda)}\right|^2\,d\lambda<\infty,
\,\frac{F}{E},\frac{F^\#}{E}\in H_2\right\}.
\end{equation*}
The space $B(E)$ with the scalar product
\begin{equation*}
[F,G]_{B(E)}=\frac{1}{\pi}\int_{\mathbb{R}}{ F(\lambda)}
\overline{G(\lambda)}\frac{d\lambda}{|E(\lambda)|^2}
\end{equation*}
is a Hilbert space. For any $z\in \mathbb{C}$ the
\emph{reproducing kernel} is introduced by the relation
\begin{equation}
\label{repr_ker} J_z(\xi):=\frac{\overline{E(z)}E(\xi)-E(\overline
z)\overline{E(\overline \xi)}}{2i(\overline z-\xi)}.
\end{equation}
Then
\begin{equation*}
F(z)=[J_z,F]_{B(E)}=\frac{1}{\pi}\int_{\mathbb{R}}{J_z(\lambda)}
\overline{F(\lambda)}\frac{d\lambda}{|E(\lambda)|^2}.
\end{equation*}
We observe that a Hermite-Biehler function $E(\lambda)$ defines
$J_z$ by (\ref{repr_ker}). The converse is also true
\cite{DMcK,DBr}:
\begin{theorem}
\label{TeorDB} Let $X$ be a Hilbert space of entire functions with
reproducing kernel such that
\begin{itemize}
\item[1)] For any $\omega\in \mathbb{C}$ the point evaluation is a
bounded functional, i.e. $|f(\omega)|\leqslant C_{\omega}\|f\|_X$,

\item[2)] if $f\in X$ then $f^\#\in X$ and $\|f\|_X=\|f^\#\|_X$,

\item[3)] if $f\in X$ and $\omega\in \mathbb{C}$ such that
$f(\omega)=0$, then $\frac{z-\overline{\omega}}{z-\omega}f(z)\in
X$ and
$\left\|\frac{z-\overline{\omega}}{z-\omega}f(z)\right\|_{X}=\|f\|_{X}$,
\end{itemize}
then $X$ is a de Branges space based on the function
\begin{equation*}
E(z)=\sqrt{\pi}(1-iz)J_i(z)\|J_i\|_X^{-1},
\end{equation*}
where $J_z$ is a reproducing kernel.
\end{theorem}

In the space $L_{2,\,\rho}(\mathbb{R})$ we take the subspace
spanned on the first $N$ polynomials generated by (\ref{Cauchy}):
\begin{equation*}
L_N:=\operatorname{Lin}\{\phi_1(\lambda),\ldots,\phi_N(\lambda)\}.
\end{equation*}
Note that $\phi_1(\lambda),\ldots,\phi_N(\lambda)$ are mutually
orthogonal in $L_{2,\,\rho}(\mathbb{R})$, see \cite{AH}. By $P_N:
L_{2,\,\rho}(\mathbb{R})\to L_{2,\,\rho}(\mathbb{R})$ we denote
the orthogonal projector in $L_{2,\,\rho}(\mathbb{R})$ onto $L_N$
acting by the rule:
\begin{equation*}
P_Na=\sum_{k=1}^N\left(a,\phi_k\right)_{L_{2,\,\rho}(\mathbb{R})}\phi_k(\lambda),\quad
a\in L_{2,\,\rho}(\mathbb{R}).
\end{equation*}
We introduce the Fourier transformation $\mathbb{R}^N\mapsto
L_{2,\,\rho}(\mathbb{R})$ by the formula
\begin{equation*}
(Fb)(\lambda)=\sum_{k=1}^Nb_k\phi_k(\lambda),\quad
b=(b_1,\ldots,b_N)\in \mathbb{R}^N.
\end{equation*}
Note that $F$ is an unitary map between $\mathbb{R}^N$ and $L_N$,
and
\begin{equation*}
b_k=\left(Fb(\lambda),\phi_k(\lambda)\right)_{L_{2,\,\rho}(\mathbb{R})}.
\end{equation*}

In accordance with the general approach proposed in
\cite{MM6,MM7}, we consider the \emph{reachable set} of the
dynamical system ({\ref{DynSyst}}):
\begin{equation*}
U^T:=W^T\mathcal{F}^T=\left\{u^f(T)\,|\, f\in
\mathcal{F}^T\right\}.
\end{equation*}
By the Lemma \ref{LemmaCont_JM} we know that
\begin{equation*}
U^T=W^T\mathcal{F}^T_1.
\end{equation*}
Then for any $f\in \mathcal{F}^T$ we can evaluate
\begin{eqnarray*}
\left(Fu^f(T)\right)(\lambda)=\sum_{k=1}^N
\int_{\mathbb{R}}\int_0^T
S_k(T-\tau,\beta)f(\tau)\,d\tau\phi_k(\beta)\,d\rho(\beta)\phi_k(\lambda)\\
=P_N\int_0^T S(T-\tau,\cdot)f(\tau)\,d\tau.
\end{eqnarray*}
We introduce the linear manifold of Fourier images of the
reachable set:
\begin{equation*}
B_N:=FU^T=\operatorname{Lin}\{\phi_1,\ldots,\phi_N\},
\end{equation*}
thus $B_N$ is a set of polynomials with complex coefficients of
the degree not grater than $N-1$.

The metric in $B_N$ is introduced by the following rule: for
$H,G\in B_N,$ such that $H=P_N\int_0^T
S(T-\tau,\cdot)h(\tau)\,d\tau,$ $G=P_N\int_0^T
S(T-\tau,\cdot)g(\tau)\,d\tau$, where $h,g\in \mathcal{F}^T_1$ we
set
\begin{equation*}
\left(H,G\right)_{B^T}:=\left(C^Th,g\right)_{\mathcal{F}^T}.
\end{equation*}
On the other hand, for $h,g\in \mathcal{F}^T_1$ we can evaluate
using the definition of $C^T$ and Fourier transformation:
\begin{eqnarray*}
%\label{sp_mes}
\left(H,G\right)_{B^T}=(C^Th,g)_{\mathcal{F}^T}=(u^h(T),u^g(T))_{\mathcal{H}^N}=\int_{\mathbb{R}}\left(Fu^h(T)\right)(\lambda)\left(Fu^h(T)\right)(\lambda)\,d\rho(\lambda)\\
=\int_{\mathbb{R}}\left(P_N\int_0^T
S(T-\tau,\cdot)h(\tau)\,d\tau\right)(\lambda) \left(P_N\int_0^T
S(T-\tau,\cdot)g(\tau)\,d\tau\right)(\lambda)\,d\rho(\lambda)\\
=\int_{\mathbb{R}}H(\lambda) G(\lambda)\,d\rho(\lambda)
\end{eqnarray*}
We note that for the systems considered in \cite{MM6,MM7} it was a
certain option in the choosing of the measure $d\rho(\lambda)$ in
the above calculations. Due to the infinite speed of wave
propagation in (\ref{DynSyst}), we do not have this option here.

We set \emph{the special control problem} for the system
(\ref{DynSyst}): to find a control $j_z\in \mathcal{F}^T_1$ which
drives (\ref{DynSyst}) to the prescribed state
\begin{equation*}
u_k^{j_z}(T)=\overline{\phi_k(z)}, \quad k=1,\ldots,N.
\end{equation*}
 at time $t=T$. Due to Theorem \ref{Krein_Th} such a control
 exists and is unique in $\mathcal{F}^T_1$. Then for such a
 control we can evaluate:
 \begin{equation*}
\left(C^Tj_z,g\right)_{\mathcal{F}^T}=\left(u^g(T),u^{j_z}(T)\right)_{\mathcal{H}^N}=\sum_{k=1}^Nu^g_k(T)\phi_k(z)=\left(Fu^g(T)\right)(z).
 \end{equation*}
Thus for
\begin{equation*}
%\label{repr_ker1}
J_z(\lambda):=\left(Fu^{j_z}(T)\right)(\lambda)
\end{equation*}
and $G(\lambda)=\left(Fu^g(T)\right)(\lambda)$ we have that
\begin{equation*}
\left(J_z,G\right)_{B_N}=\left(C^Tj_z,g\right)_{\mathcal{F}^T}=G(z).
\end{equation*}
In other words $J_z(\lambda)$ is a reproducing kernel in $B_N$.

To show that $B_N$ is a de Branges space we use the Theorem
\ref{TeorDB}, all three conditions of which are trivially
satisfied: indeed, for $G\in B_N$ such that $G=P_N\int_0^T
S(T-\tau,\cdot)g(\tau)\,d\tau$, where $g\in \mathcal{F}^T_1$ we
can evaluate:
\begin{eqnarray*}
|G(z)|=\left|(J_z,G)_{B_N}\right|=\left|\left(C^Nj_z,g\right)_{\mathcal{F}^N}\right|\\
\leqslant
\|\left(C^T\right)^{\frac{1}{2}}j_z\|_{\mathcal{F}^N}\|\left(C^N\right)^{\frac{1}{2}}g\|_{\mathcal{F}^N}
=\|\left(C^N\right)^{\frac{1}{2}}j_z\|_{\mathcal{F}^N}\|G\|_{B_N}.
\end{eqnarray*}
Clearly $G^\#$, being a polynomial is entire and
\begin{equation*}
\|G^\#\|_{B_N}=\left(\int_{-\infty}^\infty
\overline{\overline{G(\overline{\mu})}}\,\overline{G(\overline{\mu})}\,d\rho(\mu)\right)^\frac{1}{2}=\|G\|_{B_N}.
\end{equation*}
When $\omega\in \mathbb{C}$ such that $F(\omega)=0$, then
$\frac{z-\overline{\omega}}{z-\omega}F(z)$ is entire function and
\begin{equation*}
\left\|\frac{z-\overline{\omega}}{z-\omega}G(z)\right\|_{B_N}=\left(\int_{-\infty}^\infty
{\overline{\frac{z-\overline{\omega}}{z-\omega}G(z)}}\frac{z-\overline{\omega}}{z-\omega}G(z)\,d\rho(z)\right)^\frac{1}{2}=\|G\|_{B_N}.
\end{equation*}
Thus $B_N$ is a de Branges space.

\noindent{\bf Acknowledgments}

A. S. Mikhaylov and V. S. Mikhaylov were partly supported by
Volkswagen Foundation project "From Modeling and Analysis to
Approximation".

\end{document}